\documentclass[12pt,oneside,english]{amsart}
\usepackage[T1]{fontenc}
\usepackage[latin9]{inputenc}
\usepackage{amsthm}
\usepackage{amstext}
\usepackage{amssymb}
\usepackage{graphicx}

\makeatletter
\numberwithin{equation}{section}
\numberwithin{figure}{section}
\theoremstyle{plain}
\newtheorem{thm}{\protect\theoremname}
  \theoremstyle{plain}
  \newtheorem{lem}[thm]{\protect\lemmaname}
  \theoremstyle{remark}
  \newtheorem{rem}[thm]{\protect\remarkname}
  \theoremstyle{plain}
  \newtheorem{cor}[thm]{\protect\corollaryname}

\makeatother

\usepackage{babel}
  \providecommand{\corollaryname}{Corollary}
  \providecommand{\lemmaname}{Lemma}
  \providecommand{\remarkname}{Remark}
\providecommand{\theoremname}{Theorem}

\begin{document}

\title{A Mid-point Theorem for the $\bigcup$ type shape of functions}

\maketitle
\textbf{$\qquad\qquad\qquad$$\quad$Arni S.R. Srinivasa Rao} %
\footnote{\begin{center}
Part of this work was done when the author was at Center for Mathematical
Biology, Mathematical Institute, University of Oxford, England
\par\end{center}%
},%
\footnote{Comments by George E. Andrews, Pennsylvania State University and Padala
Ramu, DRDO, Pune were very helpful to improve the exposition in the
paper. Comments from anonymous reviewers are very helpul for revision.
My sincere gratitide to all of them. ASRSR is supported by funds from
the Institute of Public and Preventive Health, Georgia Regents University,
Augusta. %
}

\begin{center}
Department of Biostatistics and Epidemiology
\par\end{center}

\begin{center}
Georgia Regents University, 
\par\end{center}

\begin{center}
1120 15th Street, Augusta, Georgia, 30912. 
\par\end{center}

\begin{center}
Tel: 706-721-1295, 
\par\end{center}

\begin{center}
Fax: 706-434-7057, 
\par\end{center}

\begin{center}
Email: arrao@gru.edu
\par\end{center}

\begin{center}
\vspace{0.1cm}

\par\end{center}

\tableofcontents{}

\textbf{\large (}\textbf{\emph{\large Accepted in Bulletin of Cybernatics
and Informatics)}}{\large \par}
\begin{abstract}
A mid-point theorem is proved in an elementary way for the $\bigcup$
type shape of functions that arise out of exponential quadratic functions.
These results are inspired from epidemic patterns and growth over
a time period. \textbf{Key words:} natural numbers mapping, mean value
theorem. \textbf{MSC:} 26A09,11A99,26E99,92D30
\end{abstract}

\section{Background and motivation}

\begin{figure}
1a) \includegraphics[scale=0.52]{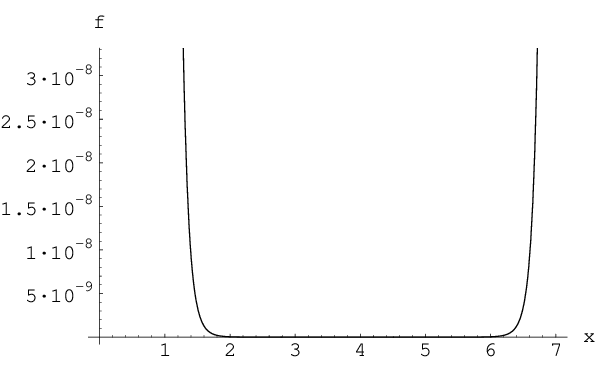} 1e)\includegraphics[scale=0.52]{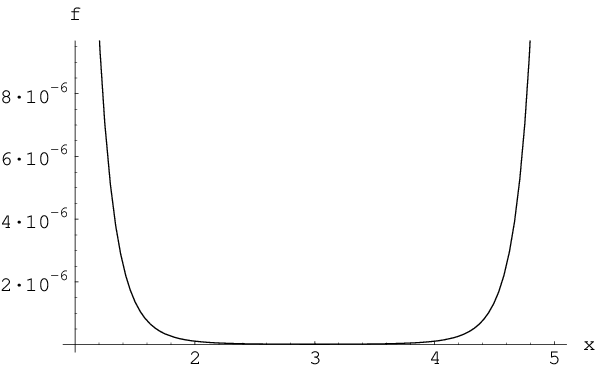}

1b) \includegraphics[scale=0.52]{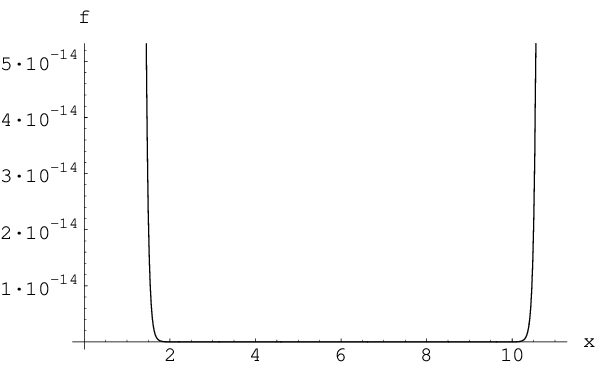} 1f) \includegraphics[scale=0.52]{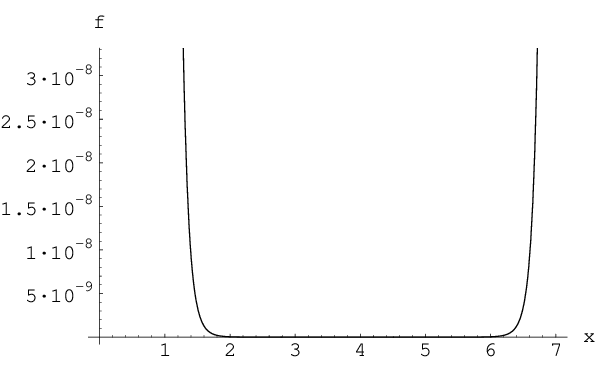}

1c) \includegraphics[scale=0.52]{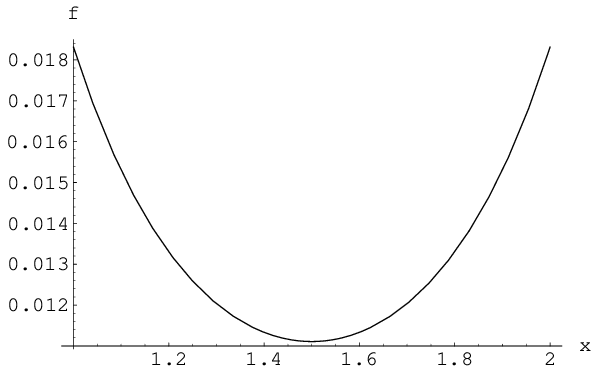} 1g) \includegraphics[scale=0.52]{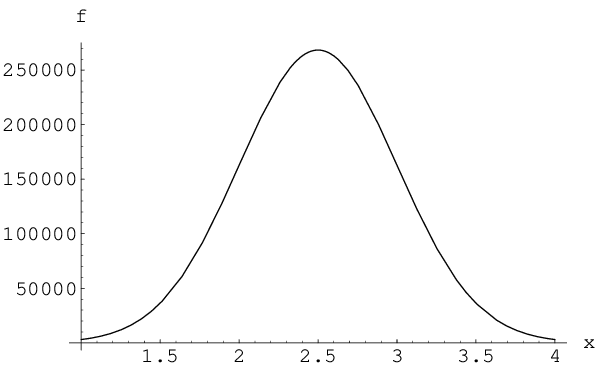}

1d) \includegraphics[scale=0.52]{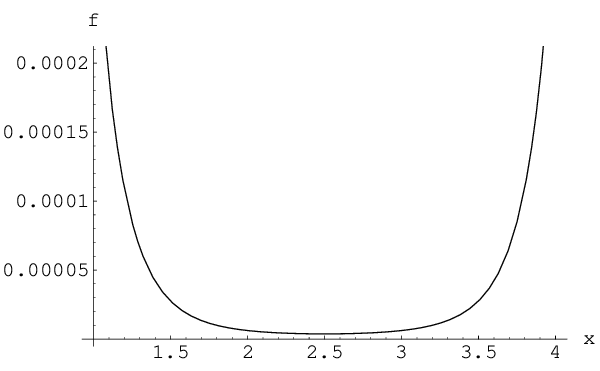} 1h) \includegraphics[scale=0.52]{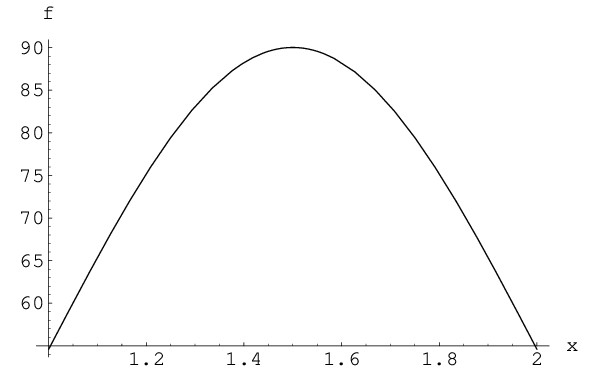}

\caption{\label{figure0}Numerical examples to demonstrate the shape of the
function $f(x)=\exp\left(ux^{2}+vx\right)$, $v=-mu^{2k-1},$ with
some $k,u,m\in\mathbb{N}$. Note that $f(1)=f(\left|A\right|)$ for
$A=\frac{u+v}{u}$. Following are combinations of $k,u,m$ in each
figure 1a) $k=2,$ $u=2,$ $m=2;$ 1b) $k=2,$ $u=2,$ $m=3;$ 1c)
$k=1,$ $u=2,$ $m=3;$ 1d) $k=1,$ $u=2,$ $m=5;$ 1e) $k=1,$ $u=2,$
$m=6;$ 1f) $k=1,$ $u=2,$ $m=8;$ 1g) $k=1,$ $u=2,$ $m=5$ (reciprocal
of the function considered); 1h) $k=1,$ $u=2,$ $m=3$ (reciprocal
of the function considered). }
\end{figure}

\begin{figure}
\includegraphics[scale=0.8]{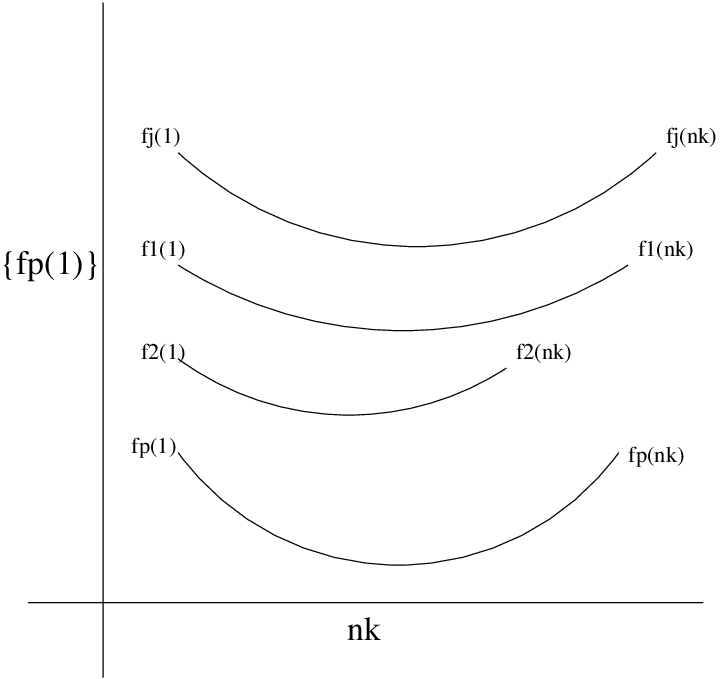}

\caption{\label{figure1}$f_{p}(1)=f_{p}(n_{k})$ for the $p^{th}$ combination
of parameters of $f$ and $n_{k}\in\mathbb{N.}$}
\end{figure}

Quadratic functions can generate variety of sizes of $\bigcup$ shaped
and $\bigcap$ shaped functions. Such kind of shapes are generally
seen, among other situations, in the growth and decay pattern of a
typical epidemic situation over certain period of years. It is often
seen, while studying the growth and decay of infections in a population
over a period of time, scientists had fitted observed epidemic data
using family of exponential or quadratic exponential functions. Quadratic
exponential functions are not only helpful in fitting the observed
data, but also often used for predicting the future course of the
epidemic \cite{key-1}. These functions consists of parameters or
constants which we estimate using the population data. In this paper,
we are not concerned various statistical methods of estimation of
parameters in the quadratic functions, but concerned in the mathematical
properties of quadratic functions, in terms of, especially in relation
to the positive integers. Some of these properties are derived while
investigating behavior of epidemic over a season in a year, over a
decade or more \cite{Rao&Kakehashi,Rao et al}. Typical epidemic data
consists of number of incidence or prevalence cases in a population
over a regular or irregular time intervals. These observations within
a given interval could either be constant or dynamic. Original work
in the direction of investigating such functions and establishing
a correspondence between natural numbers and sequence of quadratic
exponential functions in an elementary approach was inspired by realistic
situations in epidemiology\cite{key-6}. We have extended these concepts
to prove a mid-point theorem on $\bigcup$ shaped functions (see section
3). We can obtain lowest value of function under consideration between
two peaks. Suppose an epidemic pattern follows a pattern $\bigcup$,
then using this theorem we can time taken to reach lowest value of
incidence or prevalence (depending upon the context) before disease
numbers to start to grow. Further, one could try to rotate the $\bigcup$-shaped
object in a three dimensional space and obtain the volume of such
a vessel from the basic principles of Euclidean geometry. In this
paper, we have considered $\bigcup$-shaped curves and functions which
generate such curves in two dimensional space. See Figure \ref{figure0}.
Let $f=\mbox{exp}(ux^{2}+vx+w)$ be a quadratic exponential function
with (domain) $D(f)=\mathbb{N}$. For a given combination of integer
parameters (say, $\mathcal{C_{\textrm{1}}}$, for first combination
of numerical values $u$, $v$, $w$) of $f$ we establish here that
the mapping of $1\in D(f)$ will be same as mapping of some integer
$n\in D(f)$ (for the same combination of parameters). We have drawn
several curves for a combination of parameters in $f.$ If we change
the combination of parameters (say, $\mathcal{C_{\textrm{2}}}$, for
second combination of numerical values $u$, $v$, $w$), then the
resultant mapping of $1$ will be same as mapping of some $n_{1}\in D(f).$
Here $n\neq n_{1}$ and $f(1)$ for $\mathcal{C_{\textrm{1}}}$ is
not equal to $f(1)$ for $\mathcal{C_{\textrm{2}}}$ i.e. $f_{\mathcal{C_{\textrm{1}}}}(1)\neq f_{\mathcal{C_{\textrm{2}}}}(1).$
See Figure \ref{figure1} for a general idea. We construct such quadratic
exponential sequence of numbers and try to link them to the natural
numbers. The readers will also see that for some $n_{k}\in D(f)$,
the distance from $1$ to $n_{k}$ will be equal to $f(n_{k})$ in
certain conditions. We apply these facts to establish further interesting
properties of convergence and derivative of $f$. Using the principles
of mean value theorem we show that $f'(\theta)=0$ for the mid-point
$\theta\in(1,n_{k}).$ In addition to the application in epidemiology,
these results will lead to methods to compute volumes of vessels that
are of U shape in a three dimensional space. 

Consider the quadratic function, $Q_{1}(x)=u_{1}x^{2}+v_{1}x+w_{1}$
where $u_{1}>0,\, v_{1},w_{1}\in\mathbb{R}.$ Suppose $u_{1}=\frac{m_{1}}{2},\, v_{1}=\frac{n_{1}}{2}$
and $m_{1}(>0),n_{1},c_{1}\in\mathbb{N,}$ $m_{1}+n_{1}$ is even,
then $Q_{1}(x)$ is an integer for an integer $x$ \cite{key-1,key-2}.
The versatile features of quadratic function when its coefficients
are positive integers or real numbers have been popular in modeling
natural sciences \cite{key-3}. When quadratic function is taken as
an exponent to the irrational number $e$, then the resultant function
is called quadratic exponential function. Functions from such family
were widely established tools in modeling biological data \cite{key-4,key-5}.
These functions can even mimic properties of Gaussian probability
functions \cite{key-5}. Suppose $f(x)=\exp\left(ux^{2}+vx\right)$,
$v=-mu^{2k-1},$ $u,m,k\in\mathbb{N}$, then, it was proved that $f(1)=f(\left|A\right|)$
for $A=\frac{u+v}{u}$\cite{key-6}. In fact, this statement was also
proved there for $k=1$ and $u\in\mathbb{R}-\left\{ 0\right\} ,$
$m\in\mathbb{N}$ \cite{key-6}. We use these results and establish
few interesting properties of such class of exponential function.
By using Rolle's theorem we show that the derivative will be zero
at the mid-point of the interval $\left[1,\left|A\right|\right].$

Note that, $Q(1)=Q(\left|A\right|)$ when $Q(x)=ux^{2}+vx.$ Also
note $\left|A\right|=1-mu^{2k-2}$ if $(1\geq mu^{2k-2})$ or $\left|A\right|=mu^{2k-2}-1$
if $(1<mu^{2k-2}).$ We show that the above absolute value function
is necessary for deriving the main results of this paper. We begin
with some simple results.

Observe that$\left|A\right|$ is always a positive integer for any
combination of $\left\{ k,m,u\right\} \in\mathbb{N}.$ Since, $u\in\mathbb{N},\mbox{ then, }u^{2k-2}\in\mathbb{N}$
(because $m,u\in\mathbb{N}).$ Therefore $\left|A\right|\in\mathbb{N}.$ 

Additionally, whenever $u=\left(\frac{n+1}{m}\right)^{\frac{1}{2k-2}}$,
for $n\in\mathbb{N}$, then $u^{2k-2}=\frac{n+1}{m}$, which means
$\left|\frac{u+v}{u}\right|=n$. For that reason $\left|A\right|$
could be equal to every natural number for a suitable combination
of $\left\{ k,m,u\right\} ,\mbox{ where }k,m,u\in\mathbb{N}.$ For
example, if we choose $\left\{ k=1,m=2,u=2\right\} $, then, we obtain
$\left|A\right|=1$, if $\left\{ k=1,m=3,u=2\right\} $, then $\left|A\right|=2$,
if $\left\{ k=1,m=4,u=2\right\} $, then $\left|A\right|=3$, if $\left\{ k=3,m=10,u=8\right\} $,
then $\left|A\right|=40959.$ We know that $f(1)=f(\left|A\right|).$
Therefore, we modify the previous result and state that as follows: 
\begin{lem}
\label{lemma2} For any combination of $\left\{ k,m,u\right\} ,\mbox{where }k,\mbox{}m,\mbox{}u\in\mathbb{N,}$
\textup{there corresponds a} $f(1)$ such that $f(1)=f(n_{k})$ $\textrm{for some }n_{k}\in\mathbb{N}.$ \end{lem}
\begin{proof}
For every $\left|A\right|$ there corresponds a $n\in\mathbb{N}$
and $f(1)=e^{u\left(1-mu^{2k-2}\right)}$ $=f(\left|A\right|).$ Thus
$f(1)=f(n_{k})$ for some $n_{k}\in\mathbb{N}.$ 
\end{proof}
The domain of $f$ is $\mathbb{N}$. The value of $f(1)$ is not same
for every combination of $\left\{ k,m,u\right\} ,\mbox{ where }k,\mbox{}m,\mbox{}u\in\mathbb{N,}$
Readers are suggested to keep this in mind for understanding the results
presented in this work.

\section{Linking natural numbers and exponential function}
\begin{figure}
\includegraphics[scale=0.8]{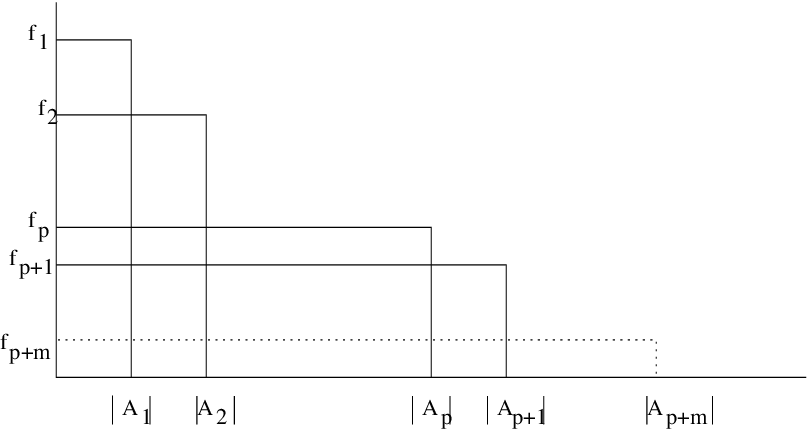}

\caption{\label{figure2}Relation between $\left|A_{p}\right|$ and $f_{p}(1).$
Vertical lines corresponding to $\left|A_{p}\right|$ are lengths
of $f_{p}(1)$ for each $p.$ }
\end{figure}

\begin{thm}
\label{theorem1}Let $N$ be even, $\left\{ m,u,k,N\right\} ,\mbox{ where }k,m,u,N\in\mathbb{N}$
and $v=\left(-1\right)^{Nk-1}mu^{Nk-1}$. Then, $f(1)=f(\left|A\right|)$.\end{thm}
\begin{proof}
$f(1)=e^{u\left(1-mu^{Nk-2}\right)}$ and $\left|A\right|=\left|\frac{u+\left(-1\right)^{Nk-1}mu^{Nk-1}}{u}\right|=mu^{Nk-2}-1$.
Therefore, $\left|A\right|\in\mathbb{N}$ for $\left\{ m,u,k,N\right\} ,\mbox{ where }k,m,u,N\in\mathbb{N}.$
Now,

\begin{eqnarray*}
f(\left|A\right|) & = & e^{u\left(mu^{Nk-2}-1\right)^{2}+v\left(mu^{Nk-2}-1\right)}\\
 & = & e^{u\left(m^{2}u^{2Nk-4}+1-2mu^{Nk-2}\right)-m^{2}u^{2Nk-3}+mu^{Nk-1}}\\
 & = & e^{u\left(1-mu^{Nk-2}\right)}
\end{eqnarray*}

Therefore, $f(1)=f(\left|A\right|).$ Since $\left|A\right|$ consists
of every element of $\mathbb{N},$ it follows that $f(1)=f(n).$\end{proof}
\begin{thm}
\label{theorem2}When $N$ is odd, $k$ is even, $v=\left(-1\right)^{Nk-1}mu^{Nk-1}$
and $\left\{ k,m,u,N\right\} ,\mbox{ where }k,\mbox{}m,\mbox{}u,\mbox{}N\in\mathbb{N}$,
then it follows that $f(1)=f(n_{k}),$ whenever $\left|A\right|\in\mathbb{N}$
and for $n_{k}\in\mathbb{N.}$\end{thm}
\begin{proof}
When $N$ is odd and $k$ is even $\left|A\right|=\left|\frac{u+\left(-1\right)^{Nk-1}mu^{Nk-1}}{u}\right|=mu^{Nk-2}-1.$
The rest of the proof can be deduced from Theorem \ref{theorem1}.\end{proof}
\begin{rem}
\label{Remark }$f{}_{1}$ denotes the $1^{st}$, $f{}_{2}$ denotes
the $2^{nd}$, and so on, $f{}_{p}$ denotes the function $f(x)=e^{\left(ux^{2}+vx\right)}$
associated with the $p^{th}$ combination of parameters $\left\{ m,u,k\right\} $,
i.e. say $\left\{ m_{1},u_{1},k_{1}\right\} ,$ $\left\{ m_{2},u_{2},k_{2}\right\} ,$...,
$\left\{ m_{p},u_{p},k_{p}\right\} ,$ then we can observe following
relations:

\begin{eqnarray*}
f_{1}(1) & = & f_{1}(2)\neq f_{1}(3)\neq...\neq f_{1}(n)\neq...\\
f_{2}(1) & \neq & f_{2}(2)=f_{2}(3)\neq...\neq f_{2}(n)\neq...\\
 & \vdots\\
\\
f_{p}(1) & = & f_{p}(2)\neq f_{p}(3)\neq...\neq f_{p}(n)=f_{p}(n+1)\neq...\\
 & \vdots
\end{eqnarray*}

Given Theorem \ref{theorem1}, suppose we denote the distance from
$\left|A\right|$ to $1$ by $\mathbb{D}$, then $f\left(\left|A\right|\right)=\mathbb{D}$,
if $m=\mathcal{Z}(1-\mathcal{Z})^{Nk-2}/\left\{ log(\mathcal{Z}-2)\right\} ^{Nk-2}.$
Here $\mathcal{Z}=mu^{Nk-2}.$ 

Suppose $\mathbb{D}_{p}$ be the logarithmic distance from $1$ to
$\left|A\right|$ for the $p^{th}$ combination of parameters, then

\begin{eqnarray*}
\mathbb{D}_{p} & = & \log\{\mathcal{Z}_{p}-2\}\\
 & = & \log\left\{ \mathcal{Z}_{p}\left(1-\frac{2}{\mathcal{Z}_{p}}\right)\right\} =\log\mathcal{Z}_{p}+\log\left(1-\frac{2}{\mathcal{Z}_{p}}\right)
\end{eqnarray*}

$\mathbb{D}_{p}$ converges for $\mathcal{Z}_{p}>2$. 

\label{remark2} The relation $f_{m}(1)=f_{m}(n)$ is unique for each
$n\in\mathbb{N}$ and $\left\{ m_{p},u_{p},k_{p}\right\} $ $\mbox{where }m_{p},u_{p},k_{p}\in\mathbb{N}.$ 

Let $N_{p+i}$ be even for $\left\{ \begin{array}{c}
p=1,2,3,...\\
i=0,1,2,...
\end{array}\right..$ If $\left|A{}_{p}\right|>\left|A_{p-1}\right|$ then $f_{p}(1)<f_{p-1}(1)$
for all $p=2,3,4,...$. This fact is demonstrated through Figure \ref{figure2}.

The number of pairs $\{f_{m}(1),\, f_{m}(n)\}$ that satisfy Remark
\ref{remark2} are countable.\end{rem}
\begin{thm}
Let $\mathcal{B}_{\sigma}(0)=\left\{ b\in\mathbb{R}^{+}:\left|b-0\right|<\sigma\right\} $
for $\sigma>0$ and $N_{p+i}$ is even for $\left\{ \begin{array}{c}
p=1,2,3,...\\
i=0,1,2,...
\end{array}\right.$, if $\left|A_{p+1}\right|>\left|A_{p}\right|$ then the sequence
$\left\{ f_{p}(1)\right\} _{p=1,2,3,...}\in\mathcal{B}_{\sigma}(0)$
for $p>M\in\mathbb{N}.$ \end{thm}
\begin{proof}
$\left|A_{p+1}\right|>\left|A_{p}\right|\Rightarrow m_{p+1}u^{N_{p+1}k_{p+1}-2}-1>$
$m_{p}u^{N_{p}k_{p}-2}-1.$

\begin{eqnarray*}
\Rightarrow\left(1-m_{p+1}u^{N_{p+1}k_{p+1}-2}\right) & < & \left(1-m_{p}u^{N_{p}k_{p}-2}\right)\\
\Rightarrow u_{p+1}\left(1-m_{p+1}u^{N_{p+1}k_{p+1}-2}\right) & < & u_{p}\left(1-m_{p}u^{N_{p}k_{p}-2}\right)\\
\Rightarrow e^{u_{p+1}\left(1-m_{p+1}u^{N_{p+1}k_{p+1}-2}\right)} & < & e^{u_{p}\left(1-m_{p}u^{N_{p}k_{p}-2}\right)}
\end{eqnarray*}

This implies $f_{p+1}(1)<f_{p}(1).$ We know that $e^{-p}\rightarrow0$
as $p\rightarrow\infty.$ Thus $f_{p}(1)\in\mathcal{B}_{\sigma}(0).$\end{proof}
\begin{cor}
Since $f_{p}(1)\in\mathcal{B}_{\sigma}(0),$ it follows from Remark
\ref{remark2} that $\left\{ f_{1}(2),\right.$ $f_{2}(3),...,f_{p}(n)$
$\left.,f_{p+1}(n+1),...\right\} $ is a convergent sequence. 
\end{cor}

\section{Mid-point Theorem}
\begin{thm}
\label{Mid-point-theorem} Suppose $\theta\in\left(1,\left|A\right|\right)$
such that $f'(\theta)=0.$ Then this $\theta$ is the mid-point of
the interval $I=\left[1,\left|A\right|\right].$\end{thm}
\begin{proof}
It is easy to verify that $f$ is continuous on $\left[1,\left|A\right|\right]$
and differentiable on $(1,\left|A\right|)$, and from Lemma \ref{lemma2}
we have $f(1)=f(\left|A\right|)$, so by Rolle's theorem there exists
a $\theta$ $\in(1,\left|A\right|)$ such that $f'(\theta)=0.$ We
have $\left|A\right|=mu^{2k-2}-1.$ Mid-point of the interval $I$
is $\frac{m}{2}u^{2k-2}.$

$f'(\theta)=e^{u\theta^{2}+v\theta}(2u\theta+v)=f(\theta)(2u\theta+v).$
Since $f'(\theta)=0,$ this means obviously $2u\theta+v=0$, because
$f(\theta)\neq0$. Therefore $\theta=\frac{-v}{2u}=\frac{m}{2}u^{2k-2}.$
Hence $\theta$ is the mid-point of the interval $I.$ For a numerical
example, consider $\left\{ k=3,\textrm{ }m=10,u=8\right\} $ as in
section 1. For this combination $\left|A\right|=40959$, mid-point
of the interval is $20480$ and $f'(20480)=0.$\end{proof}
\begin{thm}
\label{mid-point2}Suppose $\delta_{1}<\delta_{2}<...\delta_{n},$
where $\delta_{i}(>0)\in I.$ Then for a given combination of $\{k,m,u\},$
$f(1+\delta_{1})>f(1+\delta_{2})>...f(1+\delta_{n})$$=f(\left|A\right|-\delta_{n})<...<f(\left|A\right|-\delta_{2})<f(\left|A\right|-\delta_{1})$
if and only if $(1+\delta_{n})$ is a mid-point of $I$, where $\delta_{n}=\frac{m}{2}u^{2k-2}-1.$ \end{thm}
\begin{proof}
Verify easily that $f(1+\delta)=f(\left|A\right|-\delta)$ for $\delta(>0)\in I.$
Consider

\begin{eqnarray}
\qquad f(1+\delta_{1})>...>f(1+\delta_{n}) & = & f(\left|A\right|-\delta_{n})<...<f(\left|A\right|-\delta_{1})\label{one}
\end{eqnarray}

By Theorem \ref{Mid-point-theorem} we know $f'(\theta)=0$ for $\theta\in I$.
Hence $(1+\delta_{n})$ is a mid-point.

To prove converse we begin as follows. Since mid-point of the interval
$I$ is $\frac{m}{2}u^{2k-2},$ we have

\begin{eqnarray}
f(\frac{m}{2}u^{2k-2}) & = & e^{4(\frac{m}{2}u^{2k-2})^{2}-mu^{2k-1}\frac{m}{2}u^{2k-2}}\nonumber \\
 & = & e^{-\frac{m^{2}}{4}u^{4k-3}}\label{two}
\end{eqnarray}

Now for given $\delta_{n}=\frac{m}{2}u^{2k-2}-1$, we can verify that

\begin{eqnarray}
f(1+\delta_{n})=f(\left|A\right|-\delta_{n}) & = & e^{-\frac{m^{2}}{4}u^{4k-3}}\label{three}
\end{eqnarray}

Since $\delta_{1}<\delta_{2}<...\delta_{n}$ and by equations \ref{two}
and \ref{three}, the result \ref{one} is straightforward. 
\end{proof}
For large value of the distance function $\mathbb{D}$ defined, the
shape of $f$ look like the alphabet $U$ \cite{key-6}. Suppose instead
of positive integer, let $u\in\mathbb{Z^{-}}$and other parameters
$k,m$ remain as before, and if we denote resulting function as $g$,
then the shape of $g$ was shown to have mirror image of $U$ \cite{key-6}.
Based on this information and from Theorem \ref{mid-point2}, we state
the following corollary. 
\begin{cor}
Suppose $\delta_{1}<\delta_{2}<...\delta_{n},$ where $\delta_{i}(>0)\in I.$
Then for a given combination of $\{k,m,u\},$ $g(1+\delta_{1})<g(1+\delta_{2})<...g(1+\delta_{n})$$=g(\left|A\right|-\delta_{n})>...>g(\left|A\right|-\delta_{2})>g(\left|A\right|-\delta_{1})$
if and only if $(1+\delta_{n})$ is a mid-point of $I$, where $\delta_{n}=\frac{m}{2}u^{2k-2}-1.$ 
\end{cor}
Readers can also verify Darboux's theorem for $f$ on the intervals
$[1,\theta]$ and $[\theta,\left|A\right|]$ for some $f'(1)>\beta_{1}>f'(\theta)$
or $f'(\theta)<\beta_{2}<f'(\left|A\right|)$ such that $f'(\beta_{1})=\alpha_{1}$
or $f'(\beta_{2})=\alpha_{2}$ for $\alpha_{1}\in[1,\theta]$ and
$\alpha_{2}\in[\theta,\left|A\right|].$

In general, results on dynamics and periodic properties for the quadratic
function of the form $x^{2}+K$ \cite{key-7} and periodic properties
of natural numbers \cite{key-8} can be found. However, this present
note is basically deals with a correspondence between natural numbers,
quadratic exponential function, and the convergence of such functions
mapped on natural numbers constructed using $\left|A\right|.$

\begin{thebibliography}{References}
\bibitem{key-1}Pall, G. On sums of two or four values of a quadratic
function of x. Trans. Amer. Math. Soc. 34 (1932), no. 1, 98--125.

\bibitem{key-2}Pall, G. The structure of the number of representations
function in a binary quadratic form. Trans. Amer. Math. Soc. 35 (1933),
no. 2, 491--509.

\bibitem{key-3}Ojha, V. P.; Pandey, H. A population growth model
with the marriage rate as a quadratic function of time. J. Nat. Acad.
Math. India 7 (1989), no. 2, 99--104.

\bibitem{key-4}Cox, D.R and Wermuth, N. A note on the quadratic exponential
binary distribution. Biometrika 81 (1994), no. 2, 403--408.

\bibitem{key-5}McCullagh, P. Exponential mixtures and quadratic exponential
families. Biometrika 81 (1994), no. 4, 721--729.

\bibitem{key-6}Rao, Arni S.R. Srinivasa. ``U\textquotedbl{} type
functions. Bull. Inform. Cybernet. 35 (2003), no. 1-2, 35--39.

\bibitem{Rao&Kakehashi}Rao, Arni S.R. Srinivasa and Kakehashi, M
(2005). Incubation-time distribution in back-calculation applied to
HIV/AIDS data in India. Math. Biosci. Engg. 2, 2, 263-277.

\bibitem{Rao et al}Rao, Arni S. R. Srinivasa, Chen, Maggie H. Pham,
Ba' Z., Tricco A.C., Gilca V., Duval B., Krahn M.D., Bauch C.T. Cohort
effects in dynamic models and their impact on vaccination programmes:
an example from Hepatitis A. BMC Infectious Diseases (2006): 6, 174,DOI:
10.1186/1471-2334-6-174 

\bibitem{key-7}Walde, R and Russo, P. Rational periodic points of
the quadratic function $Q_{c}(x)=x^{2}+c$. Amer. Math. Monthly 101
(1994), no. 4, pp318--331.

\bibitem{key-8}Fine, N. J. Classes of periodic sequences. Illinois
J. Math. 2 (1958) 285--302. \end{thebibliography}
\end{document}